\documentclass[12pt]{article}
\usepackage{amssymb}
\usepackage{color, amsmath,amssymb, amsfonts, amstext,amsthm, latexsym}

\usepackage{amssymb, epsfig, amssymb, latexsym}
\usepackage{amsmath}
\usepackage{graphicx}
\usepackage{longtable}

\newcommand{\s}{\sigma}

\newcommand{\om}{\omega}
\newcommand{\Om}{\Omega}
\renewcommand{\phi}{\varphi}

\newcommand{\R}{{\mathbb R}}

\newcommand{\EX}{\mathbb{E}}
\newcommand{\PX}{\mathbb{P}}

\newcommand{\F}{{\mathcal{F}}}
\newcommand{\B}{{\mathcal{B}}}

\newtheorem{theorem}{Theorem}
\newtheorem{lemma}{Lemma}
\newtheorem{remark}{Remark}

\title{Synchronization of dissipative dynamical systems driven by
non-Gaussian L\'evy noises
 \footnote{This research was partly
supported by the NSF grants 0620539 and 0731201, the Cheung Kong
Scholars Program and the K. C. Wong Education Foundation.}}

 \author{Xianming Liu$^2$,   Jinqiao Duan$^1$\footnote{Corresponding author:
 duan@iit.edu},
  Jicheng Liu$^2$ and Peter E. Kloeden$^3$\\
  \\
1. Department of Applied Mathematics\\ Illinois Institute of Technology \\
  Chicago, IL 60616, USA  \\\emph{E-mail: duan@iit.edu} \\
2. School of Mathematics and Statistics\\ Huazhong University of Science and Technology \\
  Wuhan 430074, China \\
  and \\
3. Institut f\"ur Mathematik\\ Johann Wolfgang Goethe-Universit\"at\\
D-60054, Frankfurt am Main, Germany\\
   \emph{E-mail:  kloeden@math.uni-frankfurt.de} }

\begin{document}

\date{January 9, 2009}

\maketitle

\pagestyle{plain}

\begin{abstract}
Dynamical systems driven by Gaussian noises have been considered
extensively in modeling, simulation and   theory.  However,
complex systems in engineering and science are often subject to
non-Gaussian fluctuations or uncertainties. A coupled dynamical
system under non-Gaussian L\'evy noises is considered. After
discussing cocycle property, stationary orbits and random
attractors, a synchronization phenomenon is shown to occur, when
the drift terms of the coupled system satisfy certain
dissipativity and integrability conditions. The synchronization
result implies that coupled  dynamical systems  share a dynamical
feature in some asymptotic sense.

\medskip



\end{abstract}

\section{Introduction}  \label{intro}

Synchronization of coupled dynamical systems is an unbiquitous
phenomenon that has been observed in biology, physics and other areas. It concerns  coupled dynamical systems
 that share a dynamical feature in an asymptotic sense.  A descriptive account
 of its diversity of occurrence can be found in the recent book
 \cite{stro}. Recently
 Caraballo and Kloeden  \cite{Caraballo 05, Caraballo 08} have proved that
synchronization in coupled deterministic dissipative dynamical systems
persists in the presence of various Gaussian   noises (in terms of
Brownian motion), provided that appropriate concepts of random
attractors and stochastic stationary solutions are used instead of
their deterministic counterparts.

In this paper we investigate a synchronization phenomenon for
coupled dynamical systems driven   by   non-Gaussian noises (in
terms of  L\'evy motion). We show that  couple dissipative systems exhibits
synchronization for a class of L\'evy motions.

\medskip

  Gaussian processes like Brownian motion have been widely
used to model fluctuations in   engineering and science.  The  sample paths
of a particle driven by Brownian motion are continuous in
time almost surely (i.e., no jumps), the  mean square displacement
increases linearly in time (i.e., normal diffusion), and  the
probability density function decays exponentially in space (i.e.,
light tail or exponential relaxation) \cite{Oksendal}.  But some
complex phenomena involve non-Gaussian fluctuations, with peculiar
properties such as anomalous diffusion (mean square displacement
is a nonlinear power law of time) \cite{BD90} and heavy tail
(non-exponential relaxation) \cite{Yon96}. For instance, it has
been argued that diffusion in a case of geophysical turbulence
\cite{Shlesinger} is anomalous. Loosely speaking, the diffusion
process consists of a series of ``pauses", when the particle is
trapped by a coherent structure, and ``flights" or ``jumps" or
other extreme events, when the particle moves in a jet flow.
Moreover, anomalous electrical transport properties have been
observed in some amorphous materials such as insulators,
semiconductors and polymers, where transient current is
asymptotically a power law function of time \cite{SSB91,
Herrchen}.  Finally, some paleoclimatic data \cite{Dit} indicates
heavy tail distributions and some DNA data \cite{Shlesinger} shows
long range power law decay for spatial correlation.

L\'evy motions are thought to be appropriate models for
non-Gaussian processes with jumps \cite{Sato-99}. Let us recall
that  a L\'evy motion $L(t)$, or $L_t$, is a non-Gaussian process
with independent and stationary increments, i.e., increments
$\Delta L (t, \Delta t)= L(t + \Delta t)-  L(t)$ are stationary
(therefore $\Delta L$ has no statistical dependence on $t$) and
independent for any non overlapping time lags $\Delta t$.
Moreover, its sample paths are only continuous in probability,
namely, $\PX (|L(t)-L(t_0)| \geq \delta) \to 0$ as $t\to t_0$ for
any positive $\delta$. With a suitable modification
\cite{Applebaum}, these paths may be taken as   c\`{a}dl\`{a}g,
i.e., paths are continuous on the right and have limits on the
left. This continuity is weaker than the usual continuity in time.
In fact,   a  c\`{a}dl\`{a}g function  has finite or
 at most countable discontinuities on any   time interval (see, e.g.,
 p.118, \cite{Applebaum}). This generalizes the Brownian motion $B(t)$ or $B_t$,
since $B(t)$ satisfies all these three conditions, but
\emph{additionally}, (i) almost every sample path of the Brownian
motion     is continuous in time in the usual sense, and (ii) the
increments of Brownian motion are Gaussian distributed.

\medskip

This paper is organized as follows. We first recall some basic
facts about stochastic differential equations (SDEs) driven by
L\'{e}vy noise in section 2, including a fact that the solution
mappings of such SDEs generate random dynamical systems (RDS). In
section 3, we formulate the problem of synchronization of
stochastic dynamical systems driven by L\'{e}vy noises. The main
result (Theorem \ref{sync}) and an   example are presented in
section 4.

\section{Dynamical systems driven by L\'evy noises   }  \label{back}

  Dynamical systems  driven by   non-Gaussian L\'evy motions
have attracted much attention  recently \cite{Applebaum,
Schertzer}. Under certain conditions, the SDEs driven by L\'evy
motion   generate stochastic flows \cite{ Applebaum,Kunita2004},
and also generate random dynamical systems (or cocycles)  in the
sense of Arnold \cite{Arnold}.  Recently,   exit time estimates
have been investigated by Imkeller \& Pavlyukevich
\cite{ImkellerP-06, ImkellerP-08} , and Yang \& Duan
\cite{YangDuan} for SDEs driven by L\'evy motion. This shows some
qualitatively different dynamical behaviors between SDEs driven by
Gaussian and non-Gaussian noises.

\subsection{L\'evy processes}

A  L\'evy process or motion on $\mathbb{R}^d$ is characterized by
a drift parameter $\gamma \in\mathbb{R}^d$, a   covariance $d
\times d$ matrix A and a non-negative Borel measure $\nu$, defined
on $(\R^d, \mathcal{B}(\R^d))$ and concentrated on $\R^d
\setminus\{0\}$, which satisfies
\begin{equation*}
  \int_{\R^d \setminus\{0\} } (y^2 \wedge 1) \; \nu(dy) < \infty,
\end{equation*}
or equivalently
\begin{equation*}
  \int_{\R^d \setminus\{0\} } \frac{y^2}{1+y^2}\; \nu(dy) < \infty.
\end{equation*}
This measure $\nu$ is the so called  the L\'evy jump measure of
the L\'evy process $L(t)$.  Moreover L\'evy process  L$_t$ has the
following L\'{e}vy-It\^{o} decomposition
\begin{equation} \label{decomposition}
L_t=\gamma t+ B_t +\int_{|x|<1}x\mbox{\~{N}}(t,dx)+\int_{|x|\geq1}
x N(t,dx)
\end{equation}
where $N(dt,dx)$ is Poisson random measure,
\begin{equation} \label{compensated}
\mbox{\~{N}}(dt,dx)=N(dt,dx)-\nu(dx)dt
\end{equation}
is the compensated Poisson random measure of L$_t$, and  $B_t$ is
an independent Brownian motion on $\R^d$ with covariance matrix
$A$ (see \cite{Applebaum, Sato-99, ST,YangDuan}).  We   call $(A,
\nu, \gamma)$ the \emph{generating triplet}.

\bigskip

The next useful lemma provides states some important pathwise properties of $L_t$
with two-sided time $t\in \R$. Here $|\cdot |$ denotes the usual Euclidean norm in $\R^d$.

\begin{lemma}\label{lemma1} (Pathwise boundedness and convergence)\\
 Let
$L_t$ be a two-sided L\'{e}vy motion  on $\mathbb{R}^d$ for which  $\EX
|L_1|<\infty$ and $\EX L_1= \gamma$. Then

\noindent (i)\,  $\lim_{t\rightarrow \pm \infty }
\displaystyle{\frac{1}{t}}\,  L_t =\gamma. $

\noindent (ii)\,  The integrals $\int_{-\infty }^t e^{-\lambda
(t-s)}\, dL_s(\omega)$ are pathwisely uniformly bounded in
$\lambda>0$ on finite time  intervals $[T_1,T_2]$ in $\mathbb{R}$.

\noindent (iii)\,  The integrals $\int_{T_1}^t e^{-\lambda
(t-s)}\, dL_s(\omega)$ $\rightarrow$ $0$  as $\lambda$
$\rightarrow$ $\infty$, pathwise on finite time intervals  $[T_1,
T_2]$ in $\mathbb{R}$.
\end{lemma}

\begin{proof}  (i) This convergence result comes from
 the law of large numbers,   in \cite{Sato-99}, Theorem 36.5.

\noindent (ii) Due to the continuous of function $h(t)=e^{-\lambda
t}$, on integrating by parts we obtain
$$
 \int_{-\infty }^t e^{-\lambda (t-s)}\, dL_s(\omega)=
 L_t(\omega)-\lambda  \int_{-\infty }^t e^{-\lambda
 (t-s)}L_s(\omega)\, ds.
$$
Then we use (i) to conclude (ii).

(iii) Integrating again by parts, it follows that
$$
\int_{T_1}^te^{-\lambda (t-s)}\, dL_s(\omega)=
(L_t-L_{T_1})e^{-\lambda(t-T_1)}+\lambda \int_{T_1}^te^{-\lambda
(t-s)}(L_t(\omega)-L_s(\omega))\,ds,
$$
from which the result follows.
\end{proof}

\begin{remark} The assumptions on $L_t$ in  the above lemma are satisfied
by a wide class of L\'{e}vy processes, for instance, the $\alpha
$-stable symmetric L\'{e}vy motion on $\mathbb{R}^d$ with
$1<\alpha<2$. Indeed, in this case, we have
$\int_{|x|>1}|x|\nu(dx)<\infty$, and then $E|L_1|<\infty$; see
\cite{Sato-99} Theorem 25.3.
\end{remark}

Let us introduce the canonical sample space  for L\'{e}vy processes,   the
space $\Om=D(\mathbb{R},\mathbb{R}^d)$ of     c\`{a}dl\`{a}g functions, i.e.,
continuous on the right and have limits on the
left,  defined on $\mathbb{R}$ and taking values in $\mathbb{R}^d$.

If we use the usual compact-open metric, $D(\mathbb{R},\mathbb{R}^d)$
is not separable. However, it is complete and separable when
endowed with the Skorohod metric \cite{Billingsley, Situ}, in which case we
call $D(\mathbb{R},\mathbb{R}^d)$ a Skorohod space.  The Skorohod metric on $D(\mathbb{R},\mathbb{R}^d)$
is defined as
\begin{equation*}
d(x,y):=\sum_{m=1}^{\infty}\frac{1}{2^{m}}(1\wedge
d_{m}^{\circ}(x^{m},y^{m}))\qquad for\ all \ x, y\in D
\end{equation*}
where $x^{m}(t):=g_{m}(t)x(t)$, $y^{m}(t):=g_{m}(t)y(t)$ with
\begin{equation*}
g_{m}(t):= \left\{
\begin{array}{rl}
1,& \text{if } |t|\leq m-1
\\[2ex]
m-t,& \text{if } m-1\leq |t| \leq m,
\\[2ex]
0,& \text{if } |t|\geq m
\end{array} \right.
\end{equation*}
and
\begin{equation*}
d_{m}^{\circ}(x,y):=\inf_{\lambda\in\Lambda}\left\{\sup_{-m\leq s<t\leq
m}\left|\log\frac{\lambda (t)-\lambda (s)}{t-s}\right|\vee \sup_{-m\leq t\leq
m}|x(t)-y(\lambda (t))|\right\},
\end{equation*}
where $\Lambda$ denotes the set of strictly increasing, continuous
functions from $\R$ to itself.\\

Similarly, we can define  a Skorohod space on a bounded time interval  $D([T_1,
T_2],\mathbb{R}^d)$. Then, in particular,  the metric $d_{1}^{\circ}$ is the  Skorohod metric on
 $D([-1,1],\mathbb{R}^d)$.

We recall the following compactness result  (see
\cite{Billingsley}, Page 116)  in $D([T_1, T_2],\mathbb{R}^d)$.

\begin{lemma} \label{lemmaAA}
 (Ascoli-Arzela    theorem in  $D([T_1, T_2],\mathbb{R}^d)$)\\
  For $S \subset [T_1,T_2]$, let
  $$
  w_x(S)=sup \{ |x(s)-x(t)|:s,t\in S
\}
$$ and  for $0<\delta<1$ define
$$
w'_x (\delta)= \inf_{t_i}\max _{0<i\leq r}w_x \left([t_{i-1}, t_i)\right),
$$
where the infimum is taken  over
all the finite sets $\{t_i\}$ of points satisfying $ T_1$ $=$ $t_0$
$<$ $t_1$ $<$ $\ldots$ $<$ $t_r$$ =$ $T_2$ with $t_i-t_{i-1}$ $<$ $\delta$ for $i=1,2\ldots, r$.\\

Then, a  set B has compact closure in the Skorohod space  $D([T_1,
T_2],\mathbb{R}^d)$  if and only if $\sup_{x\in B} \sup_t |x(t)|<
\infty$ and $\lim_{\delta \rightarrow 0} \sup_{x\in A} w'_x
(\delta)= 0 $.
\end{lemma}

\subsection{SDE driven by L\'evy processes}

We consider the following stochastic differential equation (SDE)  driven by L\'evy motion, which has
continuous drift and Brownian motion components, namely
\begin{equation*}
  dY(t)=b(Y(t-))dt+\sigma(Y(t-))dB_t +\int_{|x|<c}F(Y(t-),x)\mbox{\~{N}}(dt,dx)
\end{equation*}
\begin{equation}\label{eq}
+\int_{|x|\geq c}G(Y(t-),x)N(dt,dx)
\end{equation}
where $\mbox{\~{N}}(dt,dx)$ and $ N(dt,dx)$ are defined above, and
the coefficients $b, \s, F, G$ are all assumed to be measurable.
Here $F$ and $G$ may be different, while the positive parameter $c$
may be different from $1$, which allows  greater generality.

We introduce  the $d\times d$ matrix
\begin{equation} \label{a}
 a(x,y)=\sigma(x)\sigma(y)^T, \qquad x, y\in\mathbb{R}^d,
\end{equation}
and define
\begin{equation*}
 \|a(x,y)\|=\sum _{i=1}^d|a_{i,i}(x,y)|.
\end{equation*}
We make the following general assumptions for the SDE \eqref{eq}:\\

\noindent  \textbf{A.1}\, There
exits $K_1>0 $ such that, for all $y_1,y_2\in \mathbb{R}^d$,
$$|b(y_1)-b(y_2)|^2+\|a(y_1,y_1)-2a(y_1,y_2)+a(y_2,y_2)\|$$
$$+\int_{|x|<1}|F(y_1,x)-F(y_2,x)|^2\nu(dx)\leq K_1|y_1-y_2|^2.$$

\noindent \textbf{A.2}\, There exits $K_2>0 $ such that, for all $y\in \mathbb{R}^d$
$$|b(y)|^2+\|a(y,y)\|+\int_{|x|<1}|F(y,x)|^2\nu(dx)\leq
K_2|1+y|^2.$$

\noindent \textbf{A.3}\,  There exits $\delta>2$ and $K_3>0 $ such that, for all
$y_1,y_2\in \mathbb{R}^d$,
$$\int_{|x|<1}|F(y_1,x)-F(y_2,x)|^p\nu(dx)\leq K_3|y_1-y_2|^p,$$
for all $2\leq p\leq \delta$.\\

From \cite {Applebaum}, Theorem 6.23, page 304, we  have the following existence and uniqueness
result for  solutions of  such SDE driven by L\'evy motion
\begin{lemma} \label{lemmaEU}
Suppose that conditions (A1) and (A2) are satisfied and that the
mapping $y\rightarrow G(y,x)$ be continuous for all $|x|\geq c$.
Then there exists a unique global c\`{a}dl\`{a}g adapted solution
of the SDE \eqref{eq}.
\end{lemma}

 Note that  a  c\`{a}dl\`{a}g solution process  has finite or
 at most countable discontinuities on any   time interval (see, e.g.,
 p.118, \cite{Applebaum}). For more details about SDEs driven by
L\'{e}vy motions, see \cite{Fujiwara, Kunita 96, Kunita2004}. Due
to the L\'{e}vy-It\^{o} decomposition \eqref{decomposition},   the
following SDE, which we consider in the sequel,
\begin{equation*}
dY(t)=f(Y(t-))dt+g(Y(t-))dL_t
\end{equation*}
is   a special case of \eqref {eq}.

\begin{remark}
The reason to take the left limit in $Y(t-)$ in the equation
\eqref {eq}  is to make sure that the c\`{a}dl\`{a}g solution
process $Y$ is predictable and unique \cite{PZ}.  For
typographical convenience, however, we will write $Y(t)$ instead
of $Y(t-)$ for the rest of this paper. Moreover, in the case of
additive noise, i.e., if the noise intensity $g(\cdot)$ does not
depend on the state $Y$, the distinction for left limit or not is
not necessary, when we consider the integral form of the equation
\eqref {eq}, as $\int_{t_0}^T f(Y(t-))dt = \int_{t_0}^T f(Y(t))dt$
for continuous $f$. In this case $f(Y(t-))$ has only countable
discontinuous points and is thus Riemann and Lebesgue integrable.

\end{remark}


\begin{remark} \label{nonglobal}
The  above global  assumptions do not hold for SDE,  which we
consider in the sequel, with a  nonlinear dissipative drift $f$
term such as $x^T f(x)$ $\leq$ $K -l |x|^2$ for some constants $K
\geq 0$ and $l>0$. However analogous global existence and
uniqueness results also hold  in this case since the dissipativity
condition prevents explosions and hence ensures otherwise local
existence is global. See \cite{Situ} for more details.
\end{remark}

\subsection{Random dynamical systems}

Following Arnold \cite{Arnold}, a random dynamical system (RDS) on
a probability space $(\Omega,\mathcal{F}, \PX)$ consists of two
ingredients: A driving   flow $\theta_t$ on the probability space
$\Om$, i.e., $\theta_t$ is a deterministic dynamical system; and a
cocycle mapping $\phi: \mathbb{R} \times\Omega\times\mathbb{R}^d
\rightarrow \mathbb{R}^d$, namely, $\phi$ satisfies the
conditions:
$$ \phi(0, \om)= id_{\R^d}, \;\;
\phi(t+s, \om)= \phi(t, \theta_s\om) \circ \phi(s, \om),
$$
for all $\om \in \Om$ and all $s, t \in \R$. This cocycle is
required to be at least measurable from the $\s-$field $\B(\R)
\times \F \times \B(\R^d)$ to the $\s-$field $\B(\R^d)$.


For random dynamical systems  driven by L\'evy noise we take
$\Omega$ $=$ $D(\mathbb{R},\mathbb{R}^d)$
with  the Skorohod metric as  the canonical sample space and
 denote by $\mathcal{F}$ $:=$ $\B(D(\mathbb{R},\mathbb{R}^d))$ the
associated Borel $\sigma$-field. Let $\mu_L$ be the (L\'evy) probability measure on $\mathcal{F}$ which
is given by the distribution of a two-sided L\'{e}vy process with
paths in $D(\mathbb{R},\mathbb{R}^d)$.

The driving system $\theta=(\theta_{t},t\in \R) $ on  $\Omega$  is defined by the
shift
\begin{equation}\label{shift}
(\theta_{t}\omega)(s):=\omega(t+s)-\omega(t).
\end{equation}
The map $(t,\omega) \rightarrow \theta_t \omega$ is continuous,
thus measurable (\cite{Arnold} page 545), and the (L\'evy) probability measure is $\theta$-invariant,
i.e.
\begin{equation*}
\mu_{L}(\theta_{t}^{-1}(A))=\mu_{L}(A)
\end{equation*}
for all $A\in\mathcal {F}$,  see \cite{Applebaum}, page 325.

\begin{lemma} \label{lemmacocycle}
(RDS generated by SDEs driven by L\'evy motion)\\
Suppose in addition to the assumptions of Lemma \ref{lemmaEU} that
condition  (A3) is satisfied. Then there exists a unique
c\`{a}dl\`{a}g adapted solution to \eqref{eq}, and the solution
mapping defines a RDS,  which is continuous in x but
c\`{a}dl\`{a}g in time.
\end{lemma}

\begin{proof}
 Let $\Phi_{s,t}$ satisfy \eqref{eq} with initial
condition $\Phi_{s,s}(y)=y$, i.e.
\begin{eqnarray}
d\Phi_{s,t}(y)& = & b(\Phi_{s,t-}(y))\, dt+\sigma(\Phi_{s,t-}(y)))\, dB_t+\int_{|x|<c}F(\Phi_{s,t-}(y)),x)\mbox{\~{N}}(dt,dx) \nonumber
\\[2ex]  \label{eq1}
& & \qquad +\int_{|x|\geq c}G(\Phi_{s,t-}(y)),x)N(dt,dx).
\end{eqnarray}
By   Theorem 6.4.2 on page 322 and  Corollary 6.4.11 on page 327 of \cite {Applebaum},  $\Phi$ is a L\'{e}vy flow, and satisfies
\begin{equation*}
\Phi_{0, s+t}(y,\omega)=\Phi_{0,t}(\Phi_{0, s}(y),
\theta_s(\omega)).
\end{equation*}
We define $\varphi: \R \times \R^d \times \Omega \rightarrow \R^d$
by
\begin{equation}\label{cocycle}
\varphi(t, y,\omega)=\Phi_{0, t}(y,\omega).
\end{equation}
It follows that
\begin{equation*}
 \varphi(t+s,y,\omega)= \varphi(t,
 \varphi(s,y,\omega),\theta_s(\omega)).
\end{equation*}
Moreover, we note that $\varphi(t,y,\omega)$ is continuous in y,
measurable in $\omega$ and c\`{a}dl\`{a}g in t. (cf.\cite
{Applebaum}, page 336). It follows that the mapping $\varphi$ is
measurable from $\R \times \R^d \times \Omega $ to $\R^d$.

With this we only need Theorem 1.3.2 and Remark 1.3.3 in
\cite{Arnold} (pages 17-20) to complete the proof.
\end{proof}

\begin{remark}
In view of Remark  \ref{nonglobal} and analogous result holds for
SDE with a  nonlinear dissipative drift term  \cite{Applebaum}.
Note that the perfection of crude discontinuous cocycles is
considered in \cite{Scheutzow}.
\end{remark}

We say that a family $\hat{A}=\{ A(\omega),\omega\in\Omega \}$ of non-empty
measurable compact subsets $A(\omega)$ of $\mathbb{R}^d$ is
$invariant$ for a RDS $(\theta,\phi)$, if $\phi(t,\omega,
A(\omega))$ $=$ $A(\theta_t \omega)$ for all $t>0$ and that it is a
random attractor if in addition it is pathwise pullback attracting
in the sense that
$$
H_d^*\left(\phi(t,\theta_{-t}\omega,D(\theta_{-t}\omega)),A(\omega)\right)\rightarrow
0 \;  \mbox{as} \;    t\rightarrow  \infty
$$
for all suitable families (called the attracting universe) of $\hat{D}$ $=$ $\{ D(\omega),\omega\in\Omega
\}$ of non-empty measurable bounded subsets $D(\omega)$ of
$\mathbb{R}^d$, where $H_d^*$ is the Hausdorff semi-distance on
$\mathbb{R}^d$.

The following result about the existence of a random attractor may
be proved similarly as in \cite{Schmalfuss2003, Caraballo 05,
Crauel, Schmalfuss, Kloeden}.

\begin{lemma} \label{lemmaattractor} (Random attractor for c\`{a}dl\`{a}g RDS)\\
 Let ($\theta,\phi$) be an RDS on
$\Omega\times\mathbb{R}^d$ and let $\phi$ be continuous in space,
but c\`{a}dl\`{a}g in time. If there exits a family $\hat{B}$ $=$ $\{
B(\omega),\omega\in\Omega \}$ of non-empty measurable compact
subsets $B(\omega)$ of $\mathbb{R}^d$ and a
$T_{\hat{D},\omega}\geq0$ such that
$$
\phi(t,\theta_{-t}\omega,D(\theta_{-t}\omega))\subset
B(\omega),\quad \forall t\geq T_{\hat{D},\omega},
$$
for all families $\hat{D}$ $=$ $\{ D(\omega),\omega\in\Omega \}$ in a
given attracting universe, then the RDS ($\theta,\phi$) has a random
attractor $\hat{A}$ $=$ $\{ A(\omega),\omega\in\Omega \}$ with the
component subsets defined for each $\omega\in\Omega$ by
$$
 A(\omega)=\bigcap_{s>0}\overline{\bigcup_{t\geq s}\phi (t,\theta_{-t}\omega,B(\theta_{-t}\omega))}
$$
Forevermore if the random attractor consists of singleton sets,
i.e $ A(\omega) =\{X^*(\omega)\}$ for some random variable $X^*$,
then $X^*_t(\omega)= X^*(\theta_t\omega)$ is a stationary
stochastic process.
\end{lemma}

We also need the following  Gronwall's lemma  from \cite{Robinson}.

\begin{lemma}\label{Gronwall}
Let $x(t)$ satisfy the differential inequality
\begin{equation*}
    \frac{d}{dt}_+x\leq g(t)x+h(t)
\end{equation*}
where $\frac{d}{dt}_+x:=lim_{h\downarrow0^+}\frac{x(t+h)-x(t)}{h}$
is right-hand derivative of  $x$.   Then
\begin{equation*}
    x(t)\leq x(0)exp[\int_0^tg(r)dr]+\int_0^t
    exp[\int_s^tg(r)dr]h(s)ds.
\end{equation*}
\end{lemma}




\section{Dissipative  synchronization}

Suppose we have two autonomous ordinary differential equations in
$\mathbb{R}^d$,
\begin{equation} \label{de1}
\frac{dx}{dt}=f(x),\quad \frac{dy}{dt}=g(y)
\end{equation}
where the vector fields $f$ and $g$ are sufficiently regular
(e.g., differentiable) to ensure the existence and uniqueness of
\emph{local} solutions, and additionally satisfy one-side
dissipative Lipschitz conditions
\begin{equation} \label{lip}
\max\left\{ \langle x_1-x_2,f(x_1)-f(x_2)\rangle, \langle
x_1-x_2,g(x_1)-g(x_2)\rangle \right\}  \leq -l{|x_1-x_2|}^2
\end{equation}
on $\mathbb{R}^d$ for some $l>0$. These dissipative Lipschitz
conditions ensure  existence and uniqueness of \emph{global}
solutions; see Remark \ref{nonglobal} above. Each of the systems
has a unique globally asymptotically stable equilibria,
$\overline{x}$ and $\overline{y}$, respectively \cite{Kloeden}.
Then,  the coupled \emph{deterministic} dynamical system in
$\mathbb{R}^{2d}$
\begin{equation} \label{de2}
\frac{dx}{dt}=f(x)+\lambda(y-x),\qquad \frac{dy}{dt}=g(x)+\lambda(x-y)
\end{equation}
with parameter $\lambda>0$ also sastisfies a one-sided  dissipative Lipschitz
condition and, hence, also  has a unique equilibrium
$(\overline{x}^\lambda,\overline{y}^\lambda)$, which is globally
asymptotically stable  \cite{Kloeden}.   Moreover,
$(\overline{x}^\lambda,\overline{y}^\lambda)\rightarrow(\overline{z},\overline{z})$
as $\lambda\rightarrow\infty$, where $\overline{z}$ is the unique
globally asymptotically stable equilibrium of the ``averaged"
system in $\mathbb{R}^{d}$
\begin{equation} \label{de3}
\frac{dz}{dt}=\frac{1}{2} \,\left(f(z)+g(z)\right).
\end{equation}
This phenomena is known as synchronization for the coupled
deterministic system \eqref{de2}. The parameter $\lambda$ often appears naturally in the context of the problem under consideraiton. For example in control theory it is a control parameter which can be chosen by the engineer, whereas in chemical reactions in thin layers separted by a membrane it is the reciprocal of the thickness of the layers, see \cite{CCK}.

Caraballo \emph{et al.} \cite{Caraballo 05, Caraballo 08} showed that
this synchronization phenomenon persists under Gaussian Brownian
noise, provided that asymptotically stable stochastic stationary
solutions are considered rather than asymptotically stable steady
state solutions. Recall that a stationary solution  $X^*$ of a SDE system may be
characterized as a    stationary orbit of the corresponding random
dynamical system   $(\theta,\phi)$(defined by the SDE system), namely,  $
\phi(t, \omega,X^*(\omega))$ $=$ $X^*(\theta_t\omega)$.

They considered a coupled system of stochastic
differential equations (SDEs) in $\mathbb{R}^{2d}$.
\begin{equation} \label{sde2}
\left\{ \begin{aligned}
         d{X_t}=\left[f({X_t})+\lambda ({Y_t}-{X_t})\right]\, dt +\alpha\,  dB_t^1,\\[2ex]
          d{Y_t}=\left[g(Y_t)+\lambda({X_t}-{Y_t})\right] \, dt+\beta\, dB_t^2.
                          \end{aligned} \right.
 \end{equation}
where $\alpha,\beta \in \mathbb{R}^d$ are constant vectors with no
components equal to zero, $B_t^1$, $B_t^2$ are independent two-sided
scalar Brownian motions, and $f, g$ satisfy the one-side
dissipative Lipschitz conditions \eqref{lip}. This coupled system
has a unique stationary solution
$(\overline{X}_t^\lambda,\overline{Y}_t^\lambda)$, which is
pathwise globally asymptotically stable. Moreover, the  coupled
system \eqref{sde2}  is synchronized to the  ``averaged" SDE in $\mathbb{R}^{d}$
\begin{equation}
\label{sde3}
dZ_t=\frac{1}{2}\, \left[f(Z_t)+g(Z_t)\right]\, dt+\frac{1}{2} \alpha\,
dB_t^1+\frac{1}{2}\beta \,dB_t^2
\end{equation}
in the sense that
$(\overline{X}_t^\lambda,\overline{Y}_t^\lambda)\rightarrow(\overline{Z}_t^\infty,\overline{Z}_t^\infty)$
as $\lambda\rightarrow\infty$, where $\overline{Z}_t^\infty$ is
the unique pathwise globally asymptotically stable stationary
solution of \eqref{sde3}.

The aim of this paper is to investigate   synchronization
  under non-Gaussian L\'{e}vy noise. In particular, we consider a coupled SDE system
  in $\mathbb{R}^d$, driven by non-Gaussian L\'{e}vy noise  in $\mathbb{R}^{2d}$
\begin{equation} \label{sde5}
\left\{ \begin{aligned}
d{X_t}=(f({X_t})+\lambda({Y_t}-{X_t}))dt +\alpha dL_t^1,
\\[1ex]
 d{Y_t}=(g(Y_t)+\lambda({X_t}-{Y_t}))dt+\beta dL_t^2,
 \end{aligned} \right.
\end{equation}
where $\alpha,\beta,f,g$ are as above, and $L_t^1, L_t^2$ are
independent two-sided  scalar L\'{e}vy  processes satisfying
conditions in Lemma \ref{lemma1}. We assume that this coupled
system defines a random dynamical system $\phi$ (i.e.., it
satisfies the assumptions in Lemma \ref{lemmacocycle} or some  generalization of it).

In addition to the one-side   Lipschitz dissipative condition
\eqref{lip} on the functions $f$ and $g$, as in  \cite{Caraballo 05} we further assume the
following integrability condition: There exists $m_0>0$ such that
for any $m\in (0,m_0]$, and any c\`{a}dl\`{a}g function $u:
\mathbb{R}\rightarrow \mathbb{R}^d$ with sub-exponential growth it
follows
\begin{equation} \label{assump}
\int_{-\infty}^t e^{ms}|f(u(s))|^{2}ds< \infty,\quad
\int_{-\infty}^t e^{ms}|g(u(s))|^{2}ds< \infty.
\end{equation}
Without loss of generality, we   assume that the
one-sided dissipative Lipschitz constant $l\leq m_0$.\\

In the next section we will show that the coupled system \eqref{sde5} has a
unique stationary solution
$(\overline{X}_t^\lambda,\overline{Y}_t^\lambda)$ which is pathwise
globally asymptotically stable with
$(\overline{X}_t^\lambda,\overline{Y}_t^\lambda)\rightarrow(\overline{Z}_t^\infty,\overline{Z}_t^\infty)$
in Skorohod metric as $\lambda\rightarrow\infty$, pathwise on
finite time-intervals $[T_1,T_2]$, where $\overline{Z}_t^\infty$
is the unique pathwise globally asymptotically stable stationary
solution of the ``averaged" SDE in $\mathbb{R}^{d}$
\begin{equation} \label{sde6}
dZ_t=\frac{1}{2}\, \left[f(Z_t)+g(Z_t)\right]\, dt+\frac{1}{2}\alpha \,
dL_t^1+\frac{1}{2}\beta \, dL_t^2.
\end{equation}

\section{Systems driven by L\'evy noise}


 For the coupled system \eqref{sde5}, we have the follow two
 lemmas   about its stationary solutions.
\begin{lemma} (Existence of stationary solutions)\\
If the Assumption \eqref{assump} holds, and f and g satisfy the
one-side Lipschitz dissipative conditions \eqref{lip}, then the
coupled stochastic system \eqref{sde5} has a unique stationary
solution.
\end{lemma}
\begin{proof}
First, the  stationary solutions  of the Langevin equations
\begin{equation} \label{sde7}
d{X_t}=-\lambda X_{t}\, dt +\alpha\, dL_t^1, \quad d{Y_t}=-\lambda
Y_{t}\, dt+\beta\,  dL_t^2
\end{equation}
 are given by
\begin{equation} \label{sde8}
\bar{X}_t^\lambda =\alpha e^{-\lambda t}\int_{-\infty}^t e^{\lambda s}
dL_t^1,\quad  \bar{Y}_t^\lambda =\beta e^{-\lambda t}\int_{-\infty}^t e^{\lambda
s} dL_t^2
\end{equation}
 The differences of the solutions of \eqref{sde5} and
these stationary solutions satisfy a system of random ordinary
differential equations, with right-hand derivative in time:
\begin{equation} \label{sde9}
\left\{ \begin{aligned} \frac{d}{dt}_+(X_t-\bar{X}_t^\lambda)
&=f(X_t)+\lambda(Y_t-X_t)+\lambda
\bar{X}_t^\lambda,\\[1ex]
\frac{d}{dt}_+(Y_t-\bar{Y}_t^\lambda)&
=g(Y_t)+\lambda(X_t-Y_t)+\lambda \bar{Y}_t^\lambda
\end{aligned} \right.
\end{equation}

 The equations \eqref{sde9} are
equivalent to
\begin{equation} \label{sde10}
\frac{d}{dt}_+U_t^\lambda=f(X_t)+\lambda(V_t^\lambda-U_t^\lambda)+\lambda
\bar{Y}_t^\lambda,\frac{d}{dt}_+V_t^\lambda=g(Y_t)+\lambda(U_t^\lambda-V_t^\lambda)+\lambda
\bar{X}_t^\lambda
\end{equation}
where $U_t^\lambda =X_t-\bar{X}_t^\lambda$ and $ V_t^\lambda
=Y_t-\bar{Y}_t^\lambda$. Thus,
$$
\frac{1}{2}\frac{d}{dt}_+(|U_t^\lambda|^2+|V_t^\lambda|^2)=
\left(U_t^\lambda,f(U_t^\lambda+\bar{X}_t^\lambda\right)-f(\bar{X}_t^\lambda))
+\left(V_t^\lambda,g(V_t^\lambda+\bar{Y}_t^\lambda\right)-g(\bar{Y}_t^\lambda))
$$
$$
\quad + \left(U_t^\lambda,f(\bar{X}_t^\lambda)+\lambda \bar{Y}_t^\lambda´\right)+ \left(V_t^\lambda,g(\bar{Y}_t^\lambda)+\lambda
\bar{X}_t^\lambda\right)-\lambda|U_t^\lambda-V_t^\lambda|^2
$$
$$
\leq
-\frac{l}{2}\left(|U_t^\lambda|^2+|V_t^\lambda|^2\right)+\frac{2}{l}\left|f(\bar{X}_t^\lambda)+\lambda
\bar{Y}_t^\lambda\right|^2+\frac{2}{l}\left|g(\bar{Y}_t^\lambda)+\lambda
\bar{X}_t^\lambda\right|^2
$$
Hence ,by Lemma \ref {Gronwall},
\begin{eqnarray*}
 |U_t^\lambda|^2+|V_t^\lambda|^2 & \leq &
 \left(|U_{t_0}^\lambda|^2+|V_{t_0}^\lambda|^2\right)
 e^{-\frac{l}{2}(t-t_0)} \\
 & & +\frac{4e^{-\frac{lt}{2}}}{l}\int_{t_0}^t e^{\frac{ls}{2}}\left[|f(\bar{X}_t^\lambda)+\lambda
\bar{Y}_t^\lambda|^2+|g(\bar{Y}_t^\lambda)+\lambda \bar{X}_t^\lambda|^2\right]\, ds
\end{eqnarray*}
This means that $|U_t^\lambda|^2+|V_t^\lambda|^2$ is pathwise
absorbed by the family $\hat{B}_{2d}^\lambda=\{
B_{2d}^\lambda(\omega),\omega\in\Omega \}$ of closed balls in
$\mathbb{R}^{2d}$ centred on the origin and of radius
$R_\lambda (\omega)$, where the $|R_\lambda
(\omega)|^2$ is given by
$$
1+\frac{4e^{-\frac{lt}{2}}}{l}\int_{-\infty}^t
e^{\frac{ls}{2}}\left[|f(\bar{X}_t^\lambda(\theta_s\omega))+\lambda
\bar{Y}_t^\lambda(\theta_s\omega)|^2+|g(\bar{Y}_t^\lambda(\theta_s\omega))+\lambda
\bar{X}_t^\lambda(\theta_s\omega)|^2\right]\, ds
$$
Hence, by Lemma \ref {lemmaattractor}, the coupled system has a
random attractor $\hat{A}^\lambda=\{
A^\lambda(\omega),\omega\in\Omega \}$ with
$A^\lambda(\omega)\subset B_{2d}^\lambda(\omega)$.

However, the difference $(\triangle X_t,\triangle
Y_t)=(X_t^1-X_t^2,Y_t^1-Y_t^2)$ of any pair of solutions satisfies
the system of random ordinary differential equations
\begin{eqnarray*}
\frac{d}{dt}_+\triangle X_t & =& f(X_t^1)-f(X_t^2)+\lambda(\triangle
Y_t-\triangle X_t), \\ \frac{d}{dt}_+\triangle Y_t & = &
g(Y_t^1)-g(Y_t^2)-\lambda(\triangle Y_t-\triangle X_t),
\end{eqnarray*}
so
\begin{eqnarray*}
 \frac{d}{dt}_+(|\triangle X_t|^2+|\triangle Y_t|^2) & = & 2(\triangle
 X_t,f(X_t^1)-f(X_t^2))+2(\triangle
 Y_t,g(Y_t^1)-g(Y_t^2))
\\
& &  \qquad -2\lambda |\triangle X_t-\triangle Y_t|^2
\\
& \leq  & -2l(|\triangle X_t|^2+|\triangle Y_t|^2)
\end{eqnarray*}
from which we obtain
$$
|\triangle X_t|^2+|\triangle Y_t|^2\leq \left(|\triangle
X_0|^2+|\triangle Y_0|^2\right) \, e^{-2lt}
$$
which means all solutions converge pathwise to each other  as
$t\rightarrow \infty$. Thus  the random attractor consists of a
singleton set  formed by an ordered pair of stationary processes
$(
\overline{X}_t^\lambda(\omega),\overline{Y}_t^\lambda(\omega))$.
\end{proof}

 \begin{remark} Using Lemma \ref{lemma1}, it can be shown that the random compact absorbing
balls $B_{2d}^\lambda(\omega)$ are contained in the common compact ball
for $\lambda \geq 1$.
\end{remark}



\begin {lemma} (A property of stationary solutions)\\
The stationary solutions of the coupled stochastic system
\eqref{sde5} have the following asymptotic behavior:
$$
\overline{X}_t^\lambda(\omega)-\overline{Y}_t^\lambda(\omega)\rightarrow 0 \qquad \mbox{as} \,\, \, \lambda\rightarrow\infty
$$
pathwise on any bounded
time-interval $[T_{1},T_{2}]$ of $\mathbb{R}$.
\end{lemma}

\begin{proof} Since
$$
 d(\overline{X}_t^\lambda-\overline{Y}_t^\lambda)= \left(-2\lambda(\overline{X}_t^\lambda-\overline{Y}_t^\lambda)+f(\overline{X}_t^\lambda) -g(\overline{Y}_t^\lambda)\right) \, dt+ \alpha dL_t^1 -\beta dL_t^2,
$$
we have
$$
d(D_t^\lambda e^{2\lambda t})= e^{2\lambda t}\left(f(\overline{X}_t^\lambda)-g(\overline{Y}_t^\lambda)\right) + \alpha
e^{2\lambda t} \,dL_t^1 -\beta e^{2\lambda t}\,dL_t^2,
$$
for with $D_t^\lambda = \overline{X}_t^\lambda-\overline{Y}_t^\lambda$, so  pathwise
\begin{eqnarray*}
|D_t^\lambda| & \leq &  e^{-2\lambda(t-T_1)}|D_{T_1}^\lambda|
+\int_{T_1}^t
e^{-2\lambda(t-s)}\left(|f(\overline{X}_s^\lambda)|+|g(\overline{Y}_s^\lambda)|\right)\, ds
\\
& & +|\alpha|\left|\int_{T_1}^t e^{-2\lambda(t-s)}\, dL_t^1\right|
+|\beta| \left|\int_{T_1}^t e^{-2\lambda(t-s)}\, dL_t^2\right|.
\end{eqnarray*}
By Lemma \ref{lemma1} we see that the radius $R_\lambda(\theta_t
\omega)$ is pathwise uniformly bounded on each bounded
time-interval $[T_1,T_2]$, so we see that the right hand of above
inequality converge to 0 as $\lambda\rightarrow \infty$ pathwise
on the bounded time-interval $[T_1,T_2]$.
\end{proof}

\bigskip

We now present the main result of this paper.

\begin{theorem} \label{sync} (Synchronization under non-Gaussian L\'evy noise)\\
Suppose  that the coupled stochastic system in $\mathbb{R}^{2d}$
\begin{equation} \label{coupled}
\left\{ \begin{aligned}
d{X_t}= \left(f({X_t})+\lambda ({Y_t}-{X_t})\right)\, dt +\alpha \, dL_t^1,\\
d{Y_t}= \left(g(Y_t)+\lambda({X_t}-{Y_t})\right)\, dt+\beta \, dL_t^2,
\end{aligned} \right.
\end{equation}
 defines a random dynamical system $(\theta,\phi)$. In addition,
assume that  $f$ and  $g$ satisfy the integrability condition
\eqref{assump} as well as the one-side   Lipschitz dissipative
condition \eqref{lip}.\\
Then the coupled stochastic system
\eqref{coupled} is synchronized to a single averaged SDE in
$\mathbb{R}^{d}$
\begin{equation} \label{sde11}
dZ_t=\frac{1}{2}\left [f(Z_t)+g(Z_t)\right]\, dt+\frac{1}{2}\alpha\,
dL_t^1+\frac{1}{2}\beta\,  dL_t^2,
\end{equation}
in the sense   that the stationary solutions of \eqref{coupled}
pathwise converge to that of \eqref{sde11}, i.e.
$(\overline{X}_t^{\lambda},\overline{Y}_t^{\lambda})$ $\rightarrow$
($Z_t^\infty, Z_t^\infty)$ in Skorohod metric on any bounded
time-interval $[T_{1},T_{2}]$  as parameter $\lambda\rightarrow
\infty$.
\end{theorem}

\begin{proof}
It is enough to demonstrate the result  for any sequence
$\lambda_n\rightarrow \infty$. Define
\begin{equation} \label{sde12}
Z_t^\lambda
:=\frac{1}{2}\left[\overline{X}_t^{\lambda}+\overline{Y}_t^{\lambda} \right],\quad   t \in \mathbb{R}.
\end{equation}
 Note that $Z_t^\lambda$ satisfies the equation
\begin{equation} \label{sde13}
dZ_t^\lambda= \frac{1}{2}\left[f(\overline{X}_t^{\lambda})+g(\overline{Y}_t^{\lambda})\right]\, dt+\frac{1}{2}\alpha
\, dL_t^1+\frac{1}{2}\beta \, dL_t^2.
\end{equation}
Also we define
$$
\overline{Z}_t:=\overline{X}_t+\overline{Y}_t,\quad  t\in
\mathbb{R},
$$
where $\overline{X}_t$ and $\overline{Y}_t$ are the (stationary)
solutions of the Langevin equations
\begin{equation} \label{sde14}
d{X_t}=-X_{t}dt +\alpha dL_t^1, \quad   d{Y_t}=-Y_{t}dt+\beta
dL_t^2,
\end{equation}
i.e.
$$
\overline{X}_t=\alpha e^{-t}\int_{-\infty}^t e^s dL_t^1, \quad
\overline{Y}_t=\beta e^{-t}\int_{-\infty}^t e^s dL_t^2.
$$
The difference $Z_t^\lambda-\overline{Z}_t$ satisfies pathwaise  a random ordinary
differential equation
$$
\frac{d}{dt}_+(Z_t^\lambda-\overline{Z}_t)=\frac{1}{2}
\left(f(\overline{X}_t^\lambda)+g(\overline{Y}_t^\lambda)\right)
+\frac{1}{2}\left(\overline{X}_t+\overline{Y}_t\right).
$$
By Lemma \ref{lemma1}, we obtain
$$
|\frac{d}{dt}_+(Z_t^\lambda(\omega)-\overline{Z}_t(\omega))|\leq
\frac{1}{2}
|f(\overline{X}_t^\lambda(\omega))+g(\overline{Y}_t^\lambda(\omega))|
+ \frac{1}{2} |\overline{X}_t(\omega)+\overline{Y}_t(\omega)|
$$
$$
\leq M_{T_{1},T_{2}}(\omega)<\infty
$$
by the c\`{a}dl\`{a}g property of the solutions and the fact that
these solutions belong to the common compact ball. We can use
Lemma \ref{lemmaAA} to conclude that for any sequence
$\lambda_n\rightarrow \infty$, there is a random subsequence
$\lambda_{n_j}(\omega)\rightarrow \infty$, such that
$Z_t^{\lambda_{n_j}}(\omega)-\overline{Z}_t(\omega)\rightarrow
Z_t^\infty(\omega)-\overline{Z}_t(\omega)$ in Skorohod metric as
$j\rightarrow\infty$. Thus $Z_t^{\lambda_{n_j}}(\omega)\rightarrow
Z_t^\infty(\omega)$ in the Skorohod metric as $j\rightarrow\infty$.
Now,
$$
Z_t^{\lambda_{n_j}}(\omega)-\overline{Y}_t^{\lambda_{n_j}}(\omega)=
\frac{\overline{X}_t^{\lambda_{n_j}}(\omega)-\overline{Y}_t^{\lambda_{n_j}}(\omega)}{2}
\rightarrow 0,
$$
$$
Z_t^{\lambda_{n_j}}(\omega)-\overline{X}_t^{\lambda_{n_j}}(\omega)=
\frac{\overline{Y}_t^{\lambda_{n_j}}(\omega)-\overline{X}_t^{\lambda_{n_j}}(\omega)}{2}
\rightarrow 0,
$$
as $ \lambda_{n_j}\rightarrow \infty $, so
$$
\overline{X}_t^{\lambda_{n_j}}(\omega)=
2Z_t^{\lambda_{n_j}}(\omega)-\overline{Y}_t^{\lambda_{n_j}}(\omega)
\rightarrow Z_t^\infty(\omega),
$$
$$
\overline{Y}_t^{\lambda_{n_j}}(\omega)=
2Z_t^{\lambda_{n_j}}(\omega)-\overline{X}_t^{\lambda_{n_j}}(\omega)
\rightarrow Z_t^\infty(\omega),
$$
as $ \lambda_{n_j}\rightarrow \infty $. Moreover,
\begin{eqnarray*}
Z_t^\lambda-\overline{Z}_t & = & Z_{T_1}^\lambda-\overline{Z}_{T_1}+
 \frac{1}{2}\int_{T_1}^t f(\overline{X}_s^{\lambda})\, ds
+ \frac{1}{2}\int_{T_1}^t g(\overline{Y}_s^{\lambda})\, ds
\\
 & & +
\frac{1}{2}\int_{T_1}^t  \overline{X}_s \,ds+ \frac{1}{2}\int_{T_1}^t
\overline{Y}_s\, ds,
\end{eqnarray*}
which converges pathwise  to
\begin{eqnarray*}
Z_t^\infty & = &  Z_{T_1}^\infty+
 \frac{1}{2}\int_{T_1}^t f(\overline{X}_s^\infty)\,ds+
 \frac{1}{2}\int_{T_1}^t g(\overline{Y}_s^\infty)\,ds
 \\ & & +
\overline{Z}_t-\overline{Z}_{T_1} +\frac{1}{2}\int_{T_1}^t
\overline{X}_s ds+ \frac{1}{2}\int_{T_1}^t \overline{Y}_s \,ds
\\
& = & Z_{T_1}^\infty+
 \frac{1}{2}\int_{T_1}^t f(\overline{X}_s^\infty)\,ds+
 \frac{1}{2}\int_{T_1}^t g(\overline{Y}_s^\infty)\, ds+
\frac{\alpha}{2}\int_{T_1}^t \,dL_s^1+
\frac{\beta}{2}\int_{T_1}^t \,dL_s^2,
\end{eqnarray*}
on the interval $[T_1,T_2]$. Therefore, $Z_t^\infty$ is a solution
of the averaged SDE \eqref{sde11} for all $t\in\mathbb{R}$. The
drift of this SDE satisfies the   dissipative one-side  condition
\eqref{lip}. It has a random attractor consisting of a singleton
set formed by a stationary   orbit, which   must be equal to
$Z_t^\infty$.

 Finally, we note that all possible subsequences of $Z_t^{\lambda_n}$
have the same pathwise limit. Thus the   full sequence
 $Z_t^{\lambda_n}$   converges to $Z_t^\infty$, as $\lambda_n\rightarrow
 \infty$. This completes the proof.
\end{proof}

\subsection {An example}

Consider two scalar  SDEs:
\begin{equation*}
    dX_t=-(X_t+1)\,dt+ dL^1_t,\quad  dY_t=-(Y_t+3)\,dt+2 \,dL^2_t,
\end{equation*}
which we   rewrite as
\begin{equation*}
    dX_t=-X_t\,dt+ dL^3_t,\quad  dY_t=-Y_t\,dt+2dL^4_t,
\end{equation*}
where $L^3_t=1+L^1_t$ and $L^4_t= 3/2+ L^2_t$.

The   corresponding coupled system  \eqref{coupled} is
 \begin{equation*}
 \left\{ \begin{aligned}
 dX_t=-X_tdt+\lambda(Y_t-X_t)\,dt+ dL^3_t,
 \\
dY_t=-Y_t\,dt+\lambda(X_t-Y_t)\, dt+2dL^4_t
\end{aligned} \right.
\end{equation*}
with the  stationary solutions
\begin{eqnarray*}
\overline{X_t^\lambda} & = & \int_{-\infty}^t e^{-(\lambda+1)(t-s)} \cosh\lambda(t-s)\,dL_s^3
+ 2 \int_{-\infty}^t e^{-(\lambda+1)(t-s)} \sinh\lambda(t-s)\,dL_s^4,
\\
\overline{Y_t^\lambda} & = & \int_{-\infty}^t e^{-(\lambda+1)(t-s)} \sinh\lambda(t-s)\,dL_s^3
+ 2 \int_{-\infty}^t e^{-(\lambda+1)(t-s)}\cosh\lambda(t-s)\,dL_s^4.
 \end{eqnarray*}
Let $\lambda\rightarrow \infty$, then
\begin{equation*}
(\overline{X_t^\lambda},\overline{Y_t^\lambda})\rightarrow
(Z_t^\infty,Z_t^\infty),
\end{equation*}
where $Z_t^\infty$, given by
\begin{equation*}
Z_t^\infty=\int_{-\infty}^t \frac{1}{2}e^{-(t-s)}\,dL_s^3
           +\int_{-\infty}^t e^{-(t-s)}\,dL_s^4,
\end {equation*}
is the stationary   solution of the following averaged SDE
\begin{equation*}
    dZ_t=-Z_t\, dt+ \frac{1}{2}\,dL^3_t+ dL^4_t,
\end{equation*}
which is equivalent to the following SDE, in terms of the original
L\'evy motions $L^1$ and $L^2$:
\begin{equation*}
    dZ_t=-(Z_t+2)\,dt+ \frac{1}{2}\,dL^1_t+ dL^2_t.
\end{equation*}

\bigskip

\noindent {\bf Acknowledgements.} We would like to thank Peter Imkeller and
Bjorn Schmalfuss for helpful discussions and comments.


\end{document}